\documentclass[12 pt]{amsart}
\usepackage{amsthm, amsfonts, amssymb, color}
\usepackage{mathrsfs}
\usepackage{amsmath}
\usepackage{amstext, amsxtra}
\usepackage{txfonts}
\usepackage[colorlinks, linkcolor=black, citecolor=blue, pagebackref, hypertexnames=false]{hyperref}
\usepackage{comment}

\allowdisplaybreaks
\usepackage{geometry}
\newtheorem{theorem}{Theorem}[section]

\newtheorem{corollary}[theorem]{Corollary}

\newtheorem{definition}[theorem]{Definition}

\newtheorem{remark}[theorem]{Remark}

\newcommand\R{\mathbb{R}}

\numberwithin{equation}{section}

\title   [] { Remark On The Notion Of  Adapted  Conformal And Other Estimates }

\author{   Avy Soffer }

\address{Avy Soffer, Rutgers University, Mathematics Department.}
\email{soffer@math.rutgers.edu}
\thanks{
A. Soffer is supported in part by NSF grant DMS-1600749 and by NSFC11671163
}

\begin{document}

\begin{abstract}
I describe a way to modify the multipliers of a-priori estimates, so as to include potential perturbations of the Laplacian.
\end {abstract}
\maketitle
\textbf{Introduction}\noindent

 A-priory estimates play a crucial role in the study of PDE. In the context of Dispersive and Hyperbolic equations this include the classical Morawetz estimate, the Conformal and Dilation identities and of course conservation Laws.
 The key to getting such estimates is a sign condition, which imply the monotonic flow in time of certain functionals of the solution. This condition fails in general for variable coefficient equations.

  For example, the Dilation identity is based on the non-negativity of the quantity $$(-i)(Au(t),Hu(t))-(-i)(Hu(t),Au(t))\geq 0,$$
 where $A:-i(x\cdot \nabla_x +\nabla_x\cdot x)/2.$  $\quad A$ is the Dilation operator.

 When $H$ is $-\Delta,$ this follows from the formal commutator identity  $ i[-\Delta,A]=2(-\Delta).$

 When the potential is present we get an extra term $i[V(x),A]=-x\cdot \nabla_x V.$ Clearly, this term has no sign in general, unless
 the potential is repulsive. See however \cite{SofX, B-S}. In some cases the sum of the potential contribution and $2(-\Delta)$ is non-negative. See e.g. \cite{B-S} in the case of one dimension.
 But that is not the case in general. Yet, the Physics community refers to any non-negative potential as repulsive!

 In recent works \cite{GLS} and \cite{LSof}, a new approach is developed to extend some of the classical positive commutator methods  to include all the continuous spectral part of $H,$ up to the threshold. This is done by constructing a multiplier(conjugate operator) which is adapted to the potential.
 The key results were that if one can construct a s.a. operator $\tilde A$, which satisfies the commutation relation
 $$
 i[H,\tilde A] =f(H), \quad f(H)\sim H \,\, \text{for} \,
     H\sim 0.
 $$
 and furthermore regularity conditions are satisfied, in particular $H$ should be of class $C^1(\tilde A)$,
 then pointwise decay estimates in time follow, for all initial data in the continuous spectral part of the Hamiltonian.

 Then it was shown that for a given $H$ for which Local-Decay (LD) estimates hold, it is possible to construct an \emph{adapted operator} $\tilde A.$

 In this note I will follow this circle of ideas to extend the Conformal and other identities to  the general case, when the Laplacian part is replaced by $-\Delta +V(x)+W(x,t).$
 Some consequences will be derived for both linear and nonlinear equations.\par

 First we introduce some notation:
 Let $H=-\Delta +V(x)$ be a Schr\"odiger operator, self-adjoint on $L^2(\R^n).$ A point in $\R^n$ is denoted by $x$, and $p \equiv -i\nabla_x,$

 $$A=(1/2)(x\cdot p+p\cdot x),\quad x\cdot p\equiv \sum_{1}^{n}x_i p_i.$$

 Then we have the following identities:
 If $\psi(t)$ solves the Schr\"odinger equation,
 \begin{align}
 &i\frac{\partial\psi}{\partial t}=H\psi\\
 &H=-\Delta+V(x)\\
 &\psi(x,0) \in L^2(\R^n).
 \end{align}
 Then, 
 \begin{align}\label{Comm}
 &i[H,A]=-2\Delta-x\cdot \nabla_x V(x)\\
&\partial_t  \langle\psi(t),B(t)\,\psi(t)\rangle= \langle\psi(t),\big [i[H,B(t)]+\frac{\partial B}{\partial t}\big ]\psi(t)\rangle\\
&\equiv  \langle\psi(t), D_H \, B(t)\, \psi(t)\rangle.
\end{align}
Here $B(t)$ is a general family of operators, and $D_H$ is called the \emph{Heisenberg Derivative}.
The above identities are initially defined in the sense of sesquilinear forms on $\mathcal{S}\times\mathcal{S}.$

$\mathcal{S}$ denotes the  Schwarz space of functions on $\R^n.$

\begin{remark}
For certain propagation estimates, an \emph{upper bound} is used; in particular to prove \emph{maximal velocity bounds}. See \cite{APSS}.
\end{remark}

\begin{definition}
A self adjoint operator $H$ is called of class $C^k(A)$, if  $(H-z)^{-1},$ the resolvent of $H,$ is strongly differentiable w.r.t. $\theta$, under the flow generated by $A$:
$$
e^{i\theta A}(H-z)^{-1}e^{i\theta A},\,\,\, z\in \text{resolvent set of}\,\, H,
$$
is  a $C^k$ function of $\theta.$
\end{definition}

\section{Propagation Estimates}
Suppose a family of self-adjoint operators $B(t)$ satisfy the following estimate:
\begin{align}\label{PROB}
&\partial_t  \langle\psi(t),B(t)\psi(t)\rangle= \langle\psi(t),C^*C\, \psi(t)\rangle+g(t)\\
&g(t)\in L^1(dt), \quad C^*C\geq 0.
\end{align}
Here $\psi(t)$ solved the Schr\"odinger equation with some given Hamiltonian $H.$
We then call the family $B(t)$ a {\bf Propagation Observable} (PROB)\cite{Sig-Sof1,HSS,Sig-Sof-LD}.\par

Upon integration over time, we obtain the bound:
\begin{align}\label{PRES}
&\int_{t_0}^{T}\|C(t)\psi(t)\|^2 dt \leq \sup_{t \in [t_0,T]} \langle\psi(t),
B(t)\,\psi(t)\rangle+C_g
&C_g =\|g(t)\|_{L^1}.
\end{align}
We call this estimate {\bf Propagation Estimate} (PRES)\cite{Sig-Sof1,HSS,Sig-Sof-LD}.

\subsection{ New A-priory Estimates}
Applying the above arguments, we see that in the case that $H=-\Delta$  we have:
$$
D_H p=D_H f(p)=0, \quad D_H x=2p, D_H (x-2pt)=0.
$$
Hence, a similar calculation with
\begin{align}
&H=-\Delta+V(x)+W(x,t)\\
&B(t)=|x-2pt|^2/t+ 4t(V+W) +B_V,
\end{align}
and using equation \ref{Comm} we get
 for the Schr\"odinger equation :

\begin{theorem}[Adapted Conformal Identity]
Let the Hamiltonian of the Schr{\"o}dinger equation be given by $H(t)=-\Delta +V(x) +W(x,t),\text{on} L^2(\R^n).$
We assume that $V,W$ are bounded and differentiable, for simplicity.
Let the initial data be in $H^1\cap <x>^{-1}L^2.$
Then, the solution at time $t$ satisfies ($H(t)\equiv H+W$)
\begin{align}
i\partial_t \langle\psi(t), [|x-2pt|^2/t+ 4t(V+W) +B_V]\psi(t)\rangle&=\\
\langle\psi(t),[-\frac{1}{t}[|x-2pt|^2/t- 4t(V+W)]+i[V,-4A]+i[W,-4A+B_V]&+4t\frac{\partial W}{\partial t}+i[H,B_V]]\psi(t)\rangle\\
=\langle\psi(t),[-\frac{1}{t^2}C(t)+[4x\cdot \nabla_x V+4V]_{-}+4x\cdot \nabla_x W&+4W+4t\frac{\partial W}{\partial t}]\psi(t)\rangle.
\end{align}
\begin{align}
&i[-\Delta +V(x),B_V]\equiv -[4x\cdot \nabla_x V+4V]_{+}\\
&C(t)\equiv |x-2pt|^2 \\
&[f(x)]_{\pm}=f(x)\quad \text{for} \pm f\geq 0, \text{otherwise}\quad 0.
\end{align}
\end{theorem}
{\bf Proof}
The proof follows from a direct  computation.
Notice that the Heisenberg derivative of the conformal factor w.r.t. the free flow is zero:($H_0=-\Delta$)
$$
D_{H_0} |x-2pt|^2=0,
$$
hence
$$
D_{H_0} \frac{1}{t}|x-2pt|^2=-\frac{1}{t^2}|x-2pt|^2.
$$
Then the Heisenberg derivative of $4t(V+W)$ is $4(V+W)+4t i[-\Delta,V+W]+4t\frac{\partial W}{\partial t},$
and this commutator term will cancel out the commutator of $V+W$ with the conformal operator, up to a left-over term, which is
$$
i[V+W,-4A]+4t\frac{\partial W}{\partial t}.
$$
Finally, we construct an operator $B_V$, adaptor,
that contributes an extra term
$$
i[-\Delta+V, B_V].
$$
This extra term is chosen (the B-trick) to satisfy
$$
i[-\Delta +V(x),B_V]\equiv -[4x\cdot \nabla_x V+4V]_{+}.
$$
See \cite{LSof}.$\blacksquare$

To use the above identity,
assume first that the potential $V$ is nonnegative, $W=0,$ and $B_V$ exists and is bounded on $L^2.$ (The physicist repulsive potential) we have:
\begin{corollary}[\bf Positive Potentials]
Under the above assumptions, assume moreover that $V$ is nonnegative, $W=0,$ and $B_V$ exists and is bounded on $L^2.$
Then,
 \emph{sharp propagation estimate} holds:
$\|\, |x-2pt|\psi(t)\|^2_{L^2}+t^2\langle\psi(t),V\psi(t)\rangle \lesssim \|\psi(0)\|_{L}^2,$

with the $L$ norm is defined as the sum of the $H^1$ norm and the $L^2(\R^n, <x>^2 d^n x).$
\end{corollary}

{\bf Proof of Corollary}

  The preliminary estimate will follow if we can show that the integral over time of the LHS of (1.1) is positive, minus the initial data, and the integral of the RHS is negative up to integrable terms.
  The LHS is positive if $V$ is positive, and the $B_V$ is positive, since it is used to cancel a positive term on the RHS.
  We postpone the details of the analysis of the adaptor term to the end.
  Then, integrating both sides from $1$ to $T$, we get that:
  \begin{align}
  &(1/t)\|\, |x-2pt|\psi(t)\|^2_{L^2}+t\langle\psi(t),V\psi(t)\rangle \lesssim \|\psi(1)\|_{L}^2,\\
 & \int_1^T dt \langle\psi(t),RHS\psi(t)\rangle=-\int_1^T dt \langle\psi(t),(1/t^2)|\ |x-2pt|^2\psi(t)\rangle\\
 &+ \int_1^T dt \langle\psi(t),[4x\cdot \nabla_x V+4V]_{-}\psi(t)\rangle.
 \end{align}
\begin{remark}
One can change the initial time to zero, by using $(t+1)^{-1}$ in the definitions of the PROB.
\end{remark}
  Note that since the RHS is a sum of negative terms, they are bounded by the initial data.
 Now, we want to iterate this estimate. We choose the new PROB without the factor$1/t.$ That is we multiply the whole PROB by $t$.
 The price we pay is that on the RHS we have a \emph{positive} term, $B_V$, coming from the time derivative hitting this term.
 So, we need to show that this term is integrable:
 $ \int_0^T dt \langle\psi(t),B_V\psi(t)\rangle \lesssim 1.$ \par
  We focus on the operator $B_V.$
 First, we construct the operator and prove it is bounded on $L^2.$
 The operator $B_V$ is defined as a solution (not unique) of the following operator equation:
 $$
 i[H,B_V]=Q(x)
 $$
 for a given $Q.$
 The formal solution for this equation can be constructed on the absolutely continuous part of the operator $H$, provided
 $Q$ is sufficiently localized, and $H$ satisfies local decay estimates on the continuous spectral part.
 In case $H=-\Delta +V, V\geq 0$, we know quite generally that the above conditions are satisfied on all of $L^2$, in three or higher dimensions. In one and two dimensions one may need extra genericity conditions.
 In this case $B_V$ exists and is bounded, and is given by
 \begin{equation}
 B_V=\int_0^{\infty} e^{iHs}(-Q)e^{-iHs}P_c(H) ds.
 \end{equation}
 \begin{theorem}
 Let $B_V=\int_0^{\infty} e^{iHs}(-Q)e^{-iHs}P_c(H) ds.$
 Assume moreover, that $H$ satisfies the following local decay estimate:
 $\int_0^{\infty} \|<x>^{-\sigma} e^{-iHt}\psi\|^2_{L^2}\leq c\|\psi\|^2_{L^2},$
 for all $\psi$ in the continuous spectral subspace of $H,$ and some $\sigma.$
 Suppose moreover that $<x>^{2\sigma}(|Q(x)|+|\Delta Q|)$ is uniformly bounded by a constant.
 Then, the integral defining $B_V$ converges to a bounded operator on $L^2.$
 \end{theorem}
 {\bf Proof}
  We note that the limit defining $B_V$ exists, since
 the integral defining the operator is a  bounded form on $L^2 \times L^2$

 \begin{align}
 &\langle\psi,B_V\phi\rangle \leq c (\int_0^{\infty} \|<x>^{-\sigma} e^{-iHt}\psi\|^2_{L^2})^{1/2}\\
 &(\int_0^{\infty}\|<x>^{-\sigma} e^{-iHt}\phi\|^2_{L^2})^{1/2}\\
 &c=\|<x>^{2\sigma}Q(x)\|_{L^{\infty}},
 \end{align}
 which holds for all $\psi, \phi \in P_c(H) L^2.$
 By a consequence of Riesz Lemma, this form defines a unique bounded operator on $L^2.$
 Next, we show the commutation relation is satisfied.
\begin{align}
& i[H,B_V]= \int_0^{\infty} e^{iHs}i[H,(-Q)]e^{-iHs}P_c(H) ds\\
&=\int_0^{\infty} \frac{d}{ds} \big [e^{iHs}(-Q)e^{-iHs}\big ]P_c(H) ds\\
&=\lim_{s\to \infty}\big [e^{iHs}(-Q)e^{-iHs}\big ]P_c(H) -(-Q)P_c(H)=QP_c(H).
\end{align}
This calculation is made in the sense of quadratic forms.
The first equality follows by the same proof as the existence of $B_V$, and the regularity assumptions on $Q.$
The limit in the last equality follows since $Q$ is bounded, and vanishes at infinity.\par 

This computation shows that in the adapted conformal identity, our choice of $Q$ is such that $B_V$ is positive.
Hence we conclude that the RHS is uniformly bounded in time.
In particular, we have for $\sigma > 3/2$ and dimension $3$ :
\begin{align}
 &\|<x>^{-\sigma} e^{-iHt}\phi\|_{L^2}^2\leq c \|e^{-iHt}\phi\|_{L^6}^2\leq \frac{1}{t^2} \|(x-2pt)\phi(t)\|^2_{L^2}\\
  & \leq c \frac{1}{t}(\|<x>\phi\|^2+\|\nabla\phi\|^2).
\end{align}
We now need the sharp decay estimate for the above bound on the $L^6$ norm.
This is proved with the help of adapted Dilation identities, following \cite{LSof}.

Given a positive potential with the above conditions on decay and regularity,
we have:
\begin{theorem}[Pointwise Weighted Decay]
Let $H=-\Delta+V$, $V\geq 0$ with
 $$<x>^{\sigma}\nabla^kV, <x>^{\sigma}(x\cdot\nabla)^kV, k=0,1,2$$
  uniformly bounded on $\R^n.$
Then, for $n\geq 3$ and $\sigma=1$, the following decay estimate holds:
$$
\|<x>^{-1}e^{-iHt}P_c(H)<x>^{-1}\|\leq \frac{c}{t}.
$$
\end{theorem}
{\bf Proof}
First, we construct the adapted dilation operator $\tilde A\equiv A+B_V,$
which satisfies the identity $i[H,\tilde A]=2H.$  See \cite{GLS,LSof}.
The construction of this operator is same as the $B_V$ used above, with a different choice of $Q.$
This operator satisfies the required regularity conditions, just like $A$; in particular, $H$ is of class $C^2(\tilde A).$
This implies that
$$
[\tilde A,e^{-itH}]=-2itHe^{-itH}, \text{on}\, D(\tilde A)\times D(\tilde A).
$$
We then get:
\begin{align}
&t\langle\psi,<x>^{-1}(H+\epsilon)^{-1/2}(H+\epsilon)e^{-itH}(H+\epsilon)^{-1/2}<x>^{-1}\phi\rangle\\
&=\langle\psi,<x>^{-1}(H+\epsilon)^{-1/2}(t\epsilon)e^{-itH}(H+\epsilon)^{-1/2}<x>^{-1}\phi\rangle\\
&-(i/2)\langle\psi,<x>^{-1}(H+\epsilon)^{-1/2}[\tilde A,e^{-itH}](H+\epsilon)^{-1/2}<x>^{-1}\phi\rangle.
\end{align}
First we show that the first term on the RHS goes to zero as $\epsilon$ goes to zero, for each fixed $t.$
For this we need the key bound:
$$
\|<x>^{-1}(H+\epsilon)^{-1/2}\|+\|(H+\epsilon)^{-1/2}<x>^{-1}\|\leq c<\infty,
$$
uniformly in $\epsilon>0.$
To prove it, we note that since $V\geq 0,$  $(H+\epsilon)\geq -\Delta+\epsilon$
and therefore by application of the spectral theorem
$(H+\epsilon)^{-1/2}\leq (-\Delta+\epsilon)^{-1/2}.$
The result now follows for $n\geq3$ , since $(-\Delta+\epsilon)^{-1/2}<x>^{-1/2}$ is bounded on $L^2.$
In one and two dimensions the same application of the Hardy Littlewood Sobolev estimate holds only for a smaller power
of the Laplacian.
This leaves us with an operator of the form $\epsilon (-\Delta+\epsilon)^{-1/4+0}\leq \epsilon^a, a>0,$ which is also uniformly bounded.
The control of the second term is more complicated. However, one can easily see that by moving $\tilde A$ through $(H+\epsilon)^{-1/2},$ we get to bound as before, an expression of the form
$$
<x>^{-1}x\cdot p(H+\epsilon)^{-1/2}
$$
and the term coming from commuting $[\tilde A, (H+\epsilon)^{-1/2}].$
This commutator, while one can formally compute using the defining commutation relation
$i[H,\tilde A]=2H,$ does not hold in general, without domain and regularity assumptions on the unbounded operators involved.
In our case, using the $C^1(\tilde A)$ property of $H,$ it follows that:
$$
i[\tilde A, (H+\epsilon)^{-1}]=(H+\epsilon)^{-1}2H(H+\epsilon)^{-1}.$$
This is the second resolvent formula.
From this we can construct the commutator with the square root of the resolvent above, by using the formula
for the square root for positive operators.
Since $i[H,A]=2H$ implies that $i[g(H),A]=2Hg'(H)$,
we conclude that
$$
i[\tilde A,(H+\epsilon)^{-1/2}]=H(H+\epsilon)^{-3/2}.
$$
This operator family is uniformly bounded on $<x>^{-1}L^2$, for $n\geq3$ as before.

\begin{remark}
In one and two dimensions it is not. So in one and two dimensions one needs to impose further conditions on the initial data,
and the spectrum at zero of $H.$
\end{remark}
Taking the limit $\epsilon \lim_{\epsilon \to 0}, \epsilon>0,$ the result follows for $n\geq3.$ $\blacksquare$

 We use this bound to estimate $B_V \psi(t).$

For any $\phi \in <x> L^2(\R^3),\quad \phi(t)\equiv e^{-iHt}\phi,$ we have
\begin{align}
&\langle\phi(t), B_V \phi(t)\rangle=\\
&\langle <x>\phi,<x>^{-1}e^{iHt}\int_0^{\infty} e^{iHs}<x>^{-2\sigma}(-Q<x>^{2\sigma})e^{-iHs} e^{-iHt}P_c(H)<x>^{-1}<x>\phi\rangle ds\\
&\leq c \|<x>\phi\|^2 \int_0^{\infty}<t+s>^{-2} ds\leq \frac{c}{t}\|<x>\phi\|^2.
\end{align}
If we now rewrite the adapted conformal identity for the PROB given by
$$
|x-2pt|^2+t^2V +tB_V,
$$
we will get a new conformal estimate as before, where now the (Integral of the) LHS is still positive.
The RHS, integrated over time, is a sum of negative terms as before, plus a positive term, coming from the derivative w.r.t. time of the $tB_V$ term on the LHS.
 This term contributes $ \int_0^{\infty}\langle \psi(t),B_V \psi(t)\rangle dt.$
 By the above estimate it is bounded by $c \ln t $ for $t$ large.
 Therefore the RHS is also bounded by  $c \ln t .    $ But,$\langle \psi(t),C(t) \psi(t)\rangle$ can not grow with time, as it would drive the LHS to $-\infty.$
 It also follows from the bound on the RHS that for large $t,$ $\langle \psi(t),t^2(-x\cdot\nabla V)_{+} \psi(t)\rangle \lesssim 1.$

The above propagation estimate is stronger than local decay estimates. The $B_V$ term is positive in this case, see below, so it also decays in time.
In particular, in 3 dimensions it implies the $L^6$ norm of the solution decays like $1/t.$
The fact that we get fast decay follows from solving the differential inequality for the derivative of the $C(t)$ part.

\subsection{General Time Independent Potentials}
The above results extend to the case of general potentials  with few more steps.
 First, we must replace
$\psi(t)$ by $P_c\psi(t)$, the projection on the continuous spectral subspace of $H.$
On the left hand side we now have a negative term coming from the negative part of $V.$
However, this is added to the conformal factor C(t), and can be shown to be positive on the range of $P_c,$
using the fast decay in space of the eigenfunctions of $H$, and assuming no bound state at zero.
First, we use the PROB based on $|x-2pt|^2/t.$
We then have
\begin{theorem}
Let $H=-\Delta +V(x)$ as before, and assume the same regularity assumptions on $V.$ Moreover, assume that $H$ has no zero eigenvalues.
Then, for $\psi(0)=P_c(H)\psi(0),$ we have
\begin{equation}
\|(x-2pt)\psi(t)\|^2/t+ \langle \psi(t),Vt \psi(t)\rangle+\int_0^{\infty} \|(x-2pt)\psi(t)\|^2/t^2 dt\leq c\|\psi(0)\|^2_L.
\end{equation}
\end{theorem}
{\bf Proof}

We construct the same PROB as before, with the adaptor operator $B_V$ used to cancel positive terms that may appear on the RHS. That is we choose $-Q(x)=4[V]_{+}+[x\cdot\nabla V]_{+}$ in the definition of $B_V$, with the domain of definition given by $P_c(H)L^2.$ As before, this implies that $B_V$ is positive.
The resulting propagation estimate is the same as before, except that now the LHS contains negative terms:
The LHS after integrating over time, is as before, given by
$$
\langle \psi(t),[4Vt+|x-2pt|^2/t+B_V] \psi(t)\rangle.
$$
Since $V$ has negative part, we need a lower bound on the LHS.
Consider then the quantity $$P_c(H)[4Vt+|x-2pt|^2/t]P_c(H)=P_c(H)e^{ix^2/4t}[4Vt+4p^2t]e^{-ix^2/4t}P_c(H).$$

Observe that on the Range of $P_c(H)$ we have that $P_c(H)(p^2+V)P_c(H)\geq0.$
Hence,
\begin{align}
&P_c(H)[4Vt+|x-2pt|^2/t]P_c(H)=P_c(H)e^{ix^2/4t}P_c(H)4(p^2+V)e^{-ix^2/4t}P_c(H)\\
&+P_c(H)(e^{ix^2/4t}-1)P_b(H)4(p^2+V)P_b(H)(e^{-ix^2/4t}-1)P_c(H)\\
&+P_c(H)(e^{ix^2/4t}-1)P_b(H)4(p^2+V)P_c(H)e^{-ix^2/4t}P_c(H)\\
&+P_c(H)e^{ix^2/4t}P_c(H)4(p^2+V)P_b(H)(e^{-ix^2/4t}-1)P_c(H),\\
&P_b\equiv I-P_c.
\end{align}
The first term on the RHS of this last equation is positive.
The other terms are higher order in powers of $1/t$, since
$$
P_b(e^{ix^2/t}-1)=P_b(ix^2/t)(1+x^2/t(1/2+\text{ higher order}))=\mathbb{O}(1)/t.
$$
We used here, (with some abuse of notation, $x^2\equiv x\cdot x$)
 that $P_b <x>^n$ is bounded, as all eigenfunctions decay exponentially, except the ones with zero energy, which we assume none exist for $H.$ \par
 
Therefore the LHS of the conformal identity at this level is positive up to terms of order 1.
The RHS at this level is negative, since the $B_V$ does not contribute.
Therefore, we get the estimates at this level, same as before. In particular, the $L^6$ norm of the solution in three dimensions decays like $1/\sqrt t.$ (in fact a bit better.)$\blacksquare$

Next, we would like to iterate the estimate, to next level. For this we need to control as before, the operator
$<x>^{-1}e^{-itH}P_c(H)<x>^{-1}$ by order $1/t$ as an operator on $L^2(\R^n).$
However, in this case, since $V$ is not positive, we can not bound $(H+\epsilon)^{-1}P_c(H)$ without further spectral assumptions.
The main new assumption we need is that
$$
P_c(H)HP_c(H)\geq P_c(H)\delta (-\Delta)P_c(H), \delta>0.
$$
In this situation we can bound $(H+\epsilon)^{-1}P_c(H)\leq \delta^{-1} P_c(H)(-\Delta+\epsilon)^{-1}/\delta P_c(H).$

This condition, is a genericity condition on the thresholds of $H$, and is expected to hold under the standard assumption of no zero energy resonances. See e.g. \cite{Sof4}, \cite{Ancona}.
With this assumption, we can complete the proof as before.
It should be noted that for the first iteration, the estimate follows without the resonance assumption.

We then have
\begin{corollary}
Let $H=-\Delta +V$ satisfy the conditions of the previous Theorem. Assume moreover that it satisfies the following genericity condition
$$
P_c(H)HP_c(H)\geq \delta  P_c(H)(-\Delta)P_c(H), \delta>0.
$$
Then the conformal estimate holds on the range of $P_c(H)<x>^{-1}L^2(\R^n).$
$n\geq 3.$
\end{corollary}

\subsection{Time Dependent Potentials}

The above methods can be extended to the case where there is an additional time dependent potential $W.$
In this case further assumptions are made on $W$, so it can be estimated in terms of the main quantities above.
The above arguments, together with further microlocalization of $B_V$ can then be used to study nonlinear dispersive equations.
The above argument applies in much the same way for the Klein-Gordon equation. It is also possible to get adapted Morawetz and Dilation identities.
The above method works in any dimension.\par
 
I will describe  a few simple cases of time dependent potentials. The source of difficulty with time dependent potentials is becoming clear, once we write the conformal identities.
Up until now two scales of the conformal identity were used:
$C(t)+4t^2V+tB_V$ and $C(t)/t+ 4tV+B_V.$ \par

Now, we will first consider the PROB \cite{Graf,Sig-Sof-INV,Lin-Strauss}
$$
C(t)/t^2 +4(V+W)+B_V/t.
$$

We assume for simplicity that $V\geq 0, W\geq 0$, and in some cases that $W$ is repulsive (but not $V$).
A direct computation of the conformal identity in this case, shows that the LHS is positive in this case, but the RHS has new terms which are not negative:
$$
\langle \psi(t),[8t^{-1}x\cdot\nabla W+4\frac{\partial W}{\partial t}+i[W,B_V]/t] \psi(t)\rangle.
$$
The first factor, coming from the gradient of $W$ is not a big problem, it may have a good sign, or controlled in some cases like before, by a fraction of the terms with good sign, or upon integration, by $\ln t.$
However the term coming from the time derivative of $W$ poses major new challenges.
 We first collect the few examples where the conformal identity closes, and imply a propagation estimate.

 {\bf Example 1- Semilinear Perturbation}

 The first example come from nonlinear dispersive equations; the case of repulsive nonlinearity
 of the form $F(|\psi|^2).$
 In this case the time derivative term is a perfect derivative:
 $$
 \langle \psi(t),\frac{\partial}{\partial t}F(|\psi|^2)\psi(t)\rangle=\frac{\partial}{\partial t}\langle \psi(t),G(|\psi|^2)\psi(t)\rangle
 $$

 for  $G=F-F_1=F-|\psi|^{-2}\int^{\psi|^2} F(z) dz.$
 If $G$ is negative, then upon integration over time of the RHS we get a negative term.
 If $G$ is positive, we move this term to the left. After integration over time, it contributes to the LHS the terms:
 $$
 \langle \psi(t),G(|\psi(t)|^2)\psi(t)\rangle
 -\langle \psi(0),G(|\psi(0)|^2)\psi(0)\rangle,
 $$
 Typically, this only modifies the original nonlinear term on the LHS.

 Next, we need to deal with the term:
 $$
 \langle \psi(t),i[W,B_V]]/t \, \psi(t)\rangle.
$$
If the initial data is small, then $W$ is small, and therefore this commutator term may be controlled by a fraction
of the conformal term. The same holds for $B_V$ small, e.g. when $V$ is small.
But there are other situations, more intricate, when one can get the smallness.
This is due to the fact, that $B_V$ is small when applied to an \emph {outgoing state}, localized away from the origin, for example.
The interesting special case is the cubic NLS in one dimension, defocusing and with a positive(but not repulsive) potential added.
In this case the RHS is:
\begin{align}\label{Cubic-nls}
&-2 \langle \psi(t),C(t)/t^3 \, \psi(t)\rangle+ \langle \psi(t),[x\cdot\partial V(x)]_{-}/t \, \psi(t)\rangle\\
& \langle \psi(t),8t^{-1}x(\partial_x |\psi(t)|^2)+4\frac{\partial |\psi(t)|^2}{\partial t} \, \psi(t)\rangle\\
& +\langle \psi(t),i[W,B_V]]/t \, \psi(t)\rangle.
\end{align}
For the derivative terms of $|\psi(t)|^2$ we use that the t derivative cancels $1/2$ of such a term on the LHS, so, it drops from the final form of the conformal estimate.
Similarly, by integration by parts, the $x$ derivative of the nonlinearity contributes a \emph{negative term} to the RHS of the form:
\begin{equation}
\langle \psi(t),-4t^{-1} |\psi(t)|^2) \psi(t)\rangle.
\end{equation}
For initial condition which is localized and small in $H^1,$ it follows that the conformal term, which now controls (combined with the $L^2$ norm) all $L^p$ norms of the solution, including $L^{\infty}$, is small for all large times.
In particular, it follows that if the norm of the operator $B_V$ is smaller than a constant of order $1$, then by bootstrap, the conformal factor goes to zero, and so is the $L^{\infty}$ norm of the solution.

This will be discussed elsewhere.
The result that follows is then, that under the above conditions on $W$ we have,
\begin{theorem}[Time-Dependent Potentials]\label{TDP1}
Let $H,W$ be as above.
Then, we have the following dispersive estimate:
\begin{align}
&\int_{1}^{\infty}\big [\|\psi(t)\|^2_{L^6} +\|\,|x-2pt|\psi(t)\|^2_{L^2}\big]\frac{dt}{t}\\
&\lesssim |\psi(1)\|^2_{L}.
\end{align}
\end{theorem}
Next we consider more general cases:

{\bf Example 2- Self similar potentials}
In this case, we consider situations in which the potential term has the property that the derivative of the potential w.r.t. time,
decays like $1/t^a.$
Clearly then, if the potential and its derivatives are localized in space, bounded, and small, and the derivative decays like $1/t$
we have that
\begin{equation}\label{SES}
 \langle \psi(t),\frac{\partial W}{\partial t} \psi(t)\rangle\leq d \langle \psi(t),t^{-1}<x>^{-\sigma} \psi(t)\rangle.
\end{equation}
If $d\ll 1$, then we can bootstrap this term. It will be controlled by the Conformal term on the RHS.
The same holds for the other $W$ terms, since they all come with a factor of $1/t.$
We therefore get that Theorem \ref{TDP1} holds.
Moreover, if $d$ is not small, the RHS terms can grow like $\ln t$ and then so is the LHS.
However, if the conformal factor $ \langle \psi(t),C(t)/t^2 \psi(t)\rangle\geq \ln t$, for $t$ large, then
The conformal term on the RHS
$$\int_{1}^{T} \langle \psi(t),t^{-2}C(t) \psi(t)\rangle \frac{dt}{t}\geq (\ln t)^2,$$
which is impossible. Therefore we conclude the following;
\begin{theorem}
Let $H=-\Delta+V(x)+W(x,t)$
be such that both $V,W$ and their derivatives up to second order(w.r.t. x) are bounded, localized in space,
and such that the derivative of $W$ w.r.t. time decays like $\delta/t, \delta\ll 1.$
Then, the conformal estimate holds, in the sense that both sides are uniformly bounded in time:

\begin{align}\label{Energy}
& \int_{1}^{T} \langle \psi(t),t^{-2}C(t) \psi(t)\rangle \frac{dt}{t} \leq c  |\psi(1)\|^2_{L}\\
&\|\psi(t)\|_{H^1}\lesssim |\psi(1)\|_{L}. \\
\end{align}
Moreover, if $\nabla_x W(x,t)=\mathbb{O}(1/t),$ the asymptotic energy projections exist:

$$
\lim_{t \to \infty} \langle \psi(t),f(H) \psi(t)\rangle=f^+(H)
$$ exists for smooth bounded functions of $H$.

If the constant $\delta$ is finite, but not small, we conclude that both sides of the conformal identity grow at most logarithmically. In particular, this implies that the conformal factor $\big \langle C(t)/t^2 \big \rangle$ is uniformly bounded.
\end{theorem}

{\bf Proof}

The first inequality was derived above.
The second inequality follows from
the energy identity
$$
\partial _t \langle \psi(t),(H+W)\psi(t)\rangle= \langle \psi(t),\frac{\partial W}{\partial t} \psi(t)\rangle\in L^1(dt).
$$

This last expression is bounded by the conformal estimate.

To prove the asymptotic energy exists, we show that the derivative is integrable:
\begin{align}
& \partial_t  \langle \psi(t),f(H) \psi(t)\rangle)=
\langle \psi(t),i[W,f(H)] \psi(t)\rangle.
\end{align}

For all $f(z)\in C^{\infty}_0$, this commutator is localized both in space and frequency.
(in fact also for functions which are constant at infinity).
Therefore, the integrability over time follows from the local decay estimates, with an extra factor of $1/t,$ since $[W,f(H)]\sim \frac{\partial W}{\partial x}G.$
As was shown before, this LD estimate follows from the control of the $L^6$ norm in three dimensions, by the conformal factor.
Similar estimates hold in other dimensions.$\blacksquare$
\begin{remark}
The fact that the asymptotic energy exists play an important role in microlocalizing,
and in scattering theory for time dependent potentials.
Typically, the assumption used \cite{Sig-Sof-LD,D-G}  is that $\nabla_x W \lesssim 1/t^{1+\epsilon}.$
\end{remark}
Up to now, we considered time dependent perturbations of special form.
Next, we consider another example, in which the conditions are rather different.
First, we notice, that the expression:
\begin{equation}\label{W}
\langle \psi(t),[8t^{-1}x\cdot\nabla W+4\frac{\partial W}{\partial t}+i[W,B_V]/t ]\psi(t)\rangle
\end{equation}
can be integrated by parts over time:
\begin{align}\label{W dot}
&\int_{1}^{T}\langle \psi(t),4\frac{\partial W}{\partial t} \psi(t)\rangle dt\\
&=\int_{1}^{T}\partial_t \langle \psi(t),4 W \psi(t)\rangle-\langle \psi(t),i[-\Delta,W] \psi(t)\rangle dt\\
&=\langle \psi(T),4 W \psi(T)\rangle-\langle \psi(1),-4 W \psi(1)\rangle-4\int_{1}^{T}\langle \psi(t),(p\cdot\nabla W+\nabla W\cdot p) \psi(t)\rangle dt.
\end{align}
The first two terms on the RHS are uniformly bounded.
The last term combined with the with the $x\cdot\nabla_x W$ of equation \ref{W} gives:
$$
\int_{1}^{T}\langle \psi(t),[(2p-x/t)\cdot\nabla W+\nabla W\cdot (2p-x/t)] \psi(t)\rangle dt.
$$
(This shows that in fact the assumption $\nabla_x W(x,t)=\mathbb{O}(1/t)$ is sufficient to control all terms.)

This integrated  expression can also be bounded by an integral over time of the conformal quantities, and we get a Gronwall's type inequality for the following expression:
$$
M(t) \leq d\int_{1}^{t}M(s) ds,
$$
$$
M(s)=\langle \psi(s),[|x-2ps|^2/s^2+<x>^{-2\sigma}]\psi(s)\rangle.
$$
Now, even if the constant $d$ is small, we only get an exponential bound on the conformal norm (and the $H^1$ norm).
The problem comes from the possibility that the derivatives grow to infinity on the support of the potential.
We therefore attempt to get some local smoothing estimates, which we can in three or more dimensions, under milder assumptions on the time dependent part.
Consider the adapted smooth version of the Morawetz estimate identity in this case:
\begin{align}\label{Morawetz}
 & \langle \psi(t),(\gamma +B_V)\psi(t)\rangle-  \langle \psi(0),(\gamma +B_V)    \psi(0)\rangle=\int_{0}^{t} \langle \psi(t),(i[H,\gamma+B_V]+i[W,\gamma]+i[W,B_V])\psi(t)\rangle,\\
&V\geq 0, i[W,\gamma]+i[W,B_V]\lesssim \delta/r^2\\
&n\geq 3.
\end{align}
  We now choose

   $$\gamma= g(r)x\cdot p+p\cdot g(r)x ,\quad [i[V(x),\gamma]]_{-} +i[-\Delta+V(x),B_V]=0,$$

$g(r)\cong <r>^{-1}$ can be chosen so that:
$$
i[-\Delta,\gamma]=-\partial_i g_{ij}\partial_j+b^2(r)\geq 0.
$$
Furthermore, $g_{ij}$ is a  positive matrix function of $r$, larger than $c/r^{1+\epsilon}, b^2(r)\geq c <r>^{-2-\epsilon}.$\par

In  three or more dimensions we therefore have:
\begin{theorem}
Suppose as before that $H=-\Delta+V(x)+W(x,t),$ on $L^2(\R^n), n\geq 3$ and $V,W$ satisfy the regularity and decay assumptions for large $x$.
Suppose moreover that the solution for some initial data in $H^{1/2}$ is uniformly bounded in $H^{1/2}.$
Then, if $\nabla_x W \lesssim t^{-a},0<a<1,$
the adapted Morawetz estimate above implies local decay and smoothing with decaying in time weight:
\begin{align}
& \int_{0}^{T} \|<x>^{-1/2-\epsilon}\nabla\psi(t)\|^2_{L^2}+\|<x>^{-1-\epsilon}\psi(t)\|^2_{L^2}dt\\
& \lesssim \sup_t\|\psi(t)\|^2_{H^{1/2}}+CT^{1-a}.
\end{align}
\end{theorem}
This estimate can then be used to prove that the conformal estimate holds with PROB $C(t)/t^{2-\theta}, \theta <1-a.$
I skip the details. It implies:
\begin{corollary}
 For any initial condition satisfying the conditions of the above  Theorem, the conformal estimate holds, and in particular we have local decay with $1/t^{\theta}$ weight, and $\|\psi(t)\|^2_{L^6}\lesssim t^{-\theta}.$
\end{corollary}


\bibliographystyle{uncrt}
\def\bibfont{\HUGE}

\end{document}